\documentclass[11pt]{article}

\newcommand{\R}{\mathbf R}
\newcommand{\Z}{\mathbf Z}

\newcommand{\sgn}{\mbox{sgn}}
\newcommand{\Pf}{\mbox{Pf}}
\newcommand{\by}{{\bf y}}

\newtheorem{Theorem}{Theorem}

\begin{document}

\title{One dimensional annihilating particle systems as extended Pfaffian point processes.}
\author{Roger Tribe, Siu Kwan Yip and Oleg Zaboronski\\Mathematics Institute, University of Warwick.\\
Coventry, CV4 7AL, UK}
\date{\today}
\maketitle
%%%%%%%%%%%%%%%%%%%%%%%%%%%%%%%%%%%%%%%%%%%%%%%%%%%%
\begin{abstract}
We prove that the multi-time particle distributions for annihilating Brownian motions, 
under the maximal entrance law on the real line, are extended Pfaffian point processes.
\end{abstract}
%%%%%%%%%%%%%%%%%%%%%%%%%%%%%%%%%%%%%%%%%%%%%%%%%%%%%
\section{Main result.}
%%%%%%%%%%%%%%%%%%%%%%%%%%%%%%%%%%%%%%%%%%%%%%%%%%%%
Consider a system of annihilating Brownian motions (ABMs) on the real line, where the particles
move independently except for instantaneous annihilation when they meet. 
Assume that the initial distribution of particles is given by a natural maximal entrance law,
which can be constructed as the infinite intensity limit of Poisson initial conditions
(see \cite{maxentrance} or \cite{ourpaper} for details). The particles,
at any fixed time $t>0$, form a simple point process on $\R$ and it is shown in \cite{ourpaper} that the 
(Lebesgue) intensities $\rho_t(z_1, z_2,\ldots, z_n)$ are given by
\[
\rho_t (z_1, z_2, \ldots, z_n) = \Pf \left[ K_t(z_i-z_j): 1\leq i,j\leq n \right]
\]
where the Pfaffian is of the $2n \times 2n$ anti-symmetric matrix constructed using 
the $2 \times 2 $ matrix kernel
\[
K_t(z) =
\left( \begin{array}{cc}
K^{11}_t(z) & K^{12}_t(z) \\
K^{21}_t(z) & K^{22}_t(z)
\end{array} \right)
= \left( \begin{array}{cc}
- t^{-1} F''(z t^{-1/2}) & - t^{-1/2} F'(z t^{-1/2}) \\
t^{-1/2} F'(zt^{-1/2}) & \sgn(z) F(|z|t^{-1/2})
\end{array} \right)
\]
and $F$ is the Gaussian error function given by
\[
F(z) = \frac{1}{2\pi^{1/2}} \int^{\infty}_z e^{-x^2/4} dx
\]
with $\sgn(z) = 1$ for $z>0$, $\sgn(z) = -1$ for $z<0$ and $\sgn(0) =0$.
Briefly, this result was derived from a Pfaffian expression for the parity interval 
probabilities (called spin variables below) for ABM's  given in \cite{bena} and \cite{ourpaper}.

This result can also expressed by saying that the positions of ABMs, at any fixed time $t> 0$,
form a Pfaffian point process on $\R$ with 
kernel $ K_t(x-y)$ (see \cite{ppp} for an introduction to Pfaffian point processes).
The purpose of this note is to show that the multi-time distributions of particles for ABMs  
can be characterised as extended Pfaffian point processes (see \cite{eppp} and references therein for
other examples). 

\vspace{.1in}

\noindent
{\bf Notation.} We write $(G_r)_{r \geq 0}$ for
the semi-group generated by convolution with the Gaussian density $g_r(z) = (2 \pi r)^{-1/2} e^{-z^2/2r}$.
\begin{Theorem} \label{T1}
Under the maximal entrance law for annihilating Brownian motions, the particle
positions at times $t>0$ form an extended Pfaffian point process, with multi-time joint intensities
\begin{equation}
\rho_{t_1 t_2 \ldots t_n}(z_1,z_2, \ldots ,z_n) = \Pf \left[ K(t_i,z_i;t_j,z_j): 1 \leq i,j \leq n \right]
\label{eppdensities}
\end{equation}
where the space-time kernel $K$ is defined as follows:  for $t>s$ and $i,j \in \{1,2\}$
\[
K^{ij}(t,x;s,y) = G_{t-s} K^{ij}_s(y-x) - 2 I_{\{i=1,j=2\}} g_{t-s}(y-x);
\]
for $t<s$  and $i \neq j \in \{1,2\}$,
\[
K^{ii}(t,x;s,y) = - K^{ii}(s,y;t,x), \quad K^{ij}(t,x;s,y) = - K^{ji}(s,y;t,x);
\]
and $K(t,x;t,y) = K_t(y-x)$. 
\end{Theorem}
\textbf{Remarks.}

\noindent 
\textbf{1.} The extra term $g_{t-s}(y-x)$ in the kernel entry $K^{12}$ is singular as $t \downarrow s$,
acting like a delta function.  This reflects the fact that that the particle at $(t,y)$ is likely to have evolved from 
the particle at $(s,x)$. Indeed consider the following heuristic approximation, for small $\epsilon$, based on the fact that
$N_s$ is a simple point measure:
\begin{eqnarray*}
E\left[N_s([z,z+\epsilon])\right] & \approx & E\left[(N_s([z,z+\epsilon]))^2\right] \\
& = & \lim_{t \downarrow s} E\left[N_s([z,z+\epsilon]) N_t([z,z+\epsilon])\right] \\
& = &  \lim_{t \downarrow s} \int^{z+\epsilon}_z \! \!    \int^{z+\epsilon}_z \!\! \rho_{st}(y,x) dx \, dy.
\end{eqnarray*}
The left hand side is $O(\epsilon)$ and it is the presence of the 
delta function that implies the same for the right hand side.

\noindent \textbf{2.} Using a thinning relation (see \cite{maxentrance} or section 2 of \cite{ourpaper}) between 
instantaneously coalescing Brownian motions (CBMs) and ABMs, it is easy to show that
under the maximal entrance law for CBMs, the particle positions, at a fixed time $t>0$, also form a Pfaffian point process 
with the kernel $K_t(\cdot)$ replaced by $2 K_t(\cdot)$.  An earlier version of this paper falsely stated that the multi-time
distributions of CBMs were also an extended Pfaffian point process. We do not have a simple description of the multi-time 
distributions for CBMs. 

\vspace{.1in}

The proof of Theorem \ref{T1} is based on the analysis of PDEs solved by certain spin variables
and particle intensities. We summarize the main steps in the proof in the next section, and refer the reader 
to \cite{ourpaper} for certain details. 

%%%%%%%%%%%%%%%%%%%%%%%%%%%%%%%%%%%%%%%%%%%%%%%%%%%

\section{Summary of the proof.}
Consider the system of ABM's on $\mathbf{R}$ under the maximal entrance law. 
We write $N_t(dx)$ for the empirical measure for the particle positions at time $t$. 
Fixing $0<r_1<r_2< \ldots<r_m$, the multi-time intensities (\ref{eppdensities}),
for $t_i \in \{r_1,\ldots,r_m\}$, act as intensities for the 
simple point process generated by $(N_{r_i}(dx):i=1,\ldots,m)$ on $m$ disjoint copies of $\R$. 
See \cite{agz} definition 4.2.3 for a careful discussion.
In particular, for almost all disjoint $z_1,\ldots,z_n$ we have
\[
\rho_{t_1 t_2 \ldots t_n}(z_1,z_2, \ldots ,z_n) = \lim_{\epsilon \downarrow 0} \epsilon^{-n}
E\left[ \prod_{i=1}^n N_{t_i}([z_i,z_i+\epsilon]) \right]. 
\]
Thus $\rho_{t_1 t_2 \ldots t_n}(z_1,z_2, \ldots ,z_n)$ acts a Lebesgue density
for the absolutely continuous part of the measure $E[ \prod_{i=1}^n N_{t_i}(dz_i)]$
(and the measure is non-singular off the set $ \cup_{i \neq j} \{z_i = z_j\}$). 
To denote such an intensity we will use the informal notation
\begin{equation} \label{intensity}
\rho_{t_1 t_2 \ldots t_n}(z_1,z_2, \ldots ,z_n) = E \left[ \prod_{i=1}^n N_{t_i}(\delta_{z_i}) \right].
\end{equation}

Let $N_t(0\!:\!x)$ be the number of particles between $0$ and $x$, for any $x \in \R$ and $t>0$. 
Borrowing the terminology from the theory of spin chains \cite{glauber} we will 
refer to the following random variables as 'spin variables':
\[
S_t(x)=(-1)^{N_t(0:x)}.
\]
In \cite{ourpaper} the multi-spin correlation function $E[S_t(x_1) \ldots S_t(x_{2m})]$ was shown to be
given by a Pfaffian. Taking derivatives in $x_1,x_3, \ldots,x_{2m-1}$ 
and then letting $x_2 \downarrow x_1, \ldots, x_{2m} \downarrow x_{2m-1}$ leads  to
the intensity $ (-2)^m \rho_t(x_1,x_3,\ldots,x_{2m-1})$. We will follow a fairly similar route, but will include an 
induction over the number of space-time points. 

The following stronger statement, giving a
Pfaffian expression for a mixed spin correlation and particle intensity, is more 
suited to an inductive proof.
\begin{Theorem} \label{T2}
Fix $\by = (y_1,\ldots,y_{2m})$ satisfying $y_1 < y_2 < \ldots < y_{2m}$ and times $0<t_1 \leq t_2 \leq \ldots \leq 
t_n \leq t$, for $m,n \geq 1$. 
Then for ABMs under the maximal entrance law, the intensities
\[
E \left[ \prod_{i=1}^n N_{t_i} (\delta_{z_i}) \; \prod_{j=1}^{2m} S_t(y_j) \right]
\]
exist, in the sense described in (\ref{intensity}), and have versions given by the $(2n+2m)$ by $(2n+2m)$ Pfaffian 
\[
%\begin{equation} \label{P2}
\Phi_{(t_i,z_i)}(t,\by) = (-2)^m 
\Pf \left[ \hat{K}(s_i,x_i;s_j,x_j): (s_i,x_i),(s_j,x_j) \in A  \right]
%\end{equation}
\]
where $A$ is the set of $2m+n$ space-time points
\[
 A =  \{(t,y_1),\ldots,(t,y_{2m})\} \cup \{(t_1,z_1), \ldots,(t_n,z_n)\} 
\] 
and the first $2m$ rows and columns correspond to the space-time points $(t,y_1), \ldots,(t,y_{2m})$. 
The kernel $\hat{K}$ is defined as follows: 
\[
\begin{array}{rclc}
\hat{K}(t_i,z_i;t_j,z_j) & = & K(t_i,z_i;t_j,z_j) &  (\mbox{$2 \times 2$ entry}), \\
\hat{K}(t,y_i;t_j,z_j) & = & \left( G_{t-t_j}  K^{21}_{t_j}(z_j-y_i), \,  G_{t-t_j} K^{22}_{t_j}(z_j-y_i)  \right) 
& (\mbox{$1 \times 2$ entry}), \\
\hat{K}(t_j,z_j;t,y_i) & = & \left( -G_{t-t_j}  K^{21}_{t_j}(z_j-y_i), \,  - G_{t-t_j} K^{22}_{t_j}(z_j-y_i)  \right)^T 
& (\mbox{$2 \times 1$ entry}), \\
\hat{K}(t,y_i;t,y_j) & = & K_t^{22}(y_j-y_i), \;\; \mbox{for $i<j$} 
& (\mbox{$1 \times 1$ entry}).
\end{array}
\]
The same result holds if one of $m$ or $n$ is zero, if we take an empty product to have value $1$.
\end{Theorem}

Note that Theorem \ref{T1} is the special case of Theorem \ref{T2} when $m=0$. Furthermore
the special case when $n=0$ was the key to the results in \cite{ourpaper}. 
Note also that the ordering of the entries corresponding to the space-time points
$(t_i,z_i)$ in the Pfaffian is not important. Indeed switching
$(t_i,z_i)$ for $(t_j,z_j)$ will switch two rows and columns at once, leaving the 
Pfaffian unchanged.

The main idea of the proof is to examine the equation solved by $\Phi_{(t_i,z_i)}(t,\by)$ 
as a function of $t \geq t_n$ and  $ \by \in \overline{V}_{2m}$, 
where $V_{2m}$ is the open cell $\{\by: y_1<y_2< \ldots < y_{2m}\}$. 

\vspace{.1in}

\noindent
{\bf Step 1.  Regularity of $\Phi$.}
The regularity of the kernel $\hat{K}$  implies that $(t,\by) \to \Phi_{(t_i,z_i)}(t,\by)$ defines a bounded function 
lying in $C^{1,2}((t_n,\infty) \times V_{2m}) \cap C((t_n,\infty) \times \overline{V}_{2m})$.
The initial condition at $t = t_n$ may have a jump discontinuity when $y_j = z_i$ for some $i,j$. However, the function $(t,\by') \to  \Phi_{(t_i,z_i)}(t,\by')$ is continuous as 
$t \downarrow t_n, \, \by' \to \by$ provided that $\{y_j\}_{j \in [1,2m]} \cap \{z_i\}_{i \in [1,n]} = \emptyset$.

\vspace{.1in}

\noindent
{\bf Step 2.  PDE for $\Phi$.}
$\Phi_{(t_i,z_i)}(t,\by) $ satisfies the heat equation 
\[
 \frac{\partial}{\partial t} \Phi_{(t_i,z_i)}(t,\by)  =  \frac12 \Delta \Phi_{(t_i,z_i)}(t,\by)  
\]
for $t > t_n$ and $\by \in V_{2m}$.
To see this consider expanding $\Phi_{(t_i,z_i)}(t,\by)$ as a finite sum arising from the terms of the Pfaffian. 
We claim each of these terms seperately solves the heat equation.
Indeed, in each summand, all occurrences of the variables $t, (y_i)$ occur inside products of terms of the form
\[
K_t^{22}(y_j-y_i), \qquad  G_{t-t_j}  K^{21}_{t_j}(z_j-y_i), \qquad G_{t-t_j} K^{22}_{t_j}(z_j-y_i).
\]
Note that each of these terms solves the heat equation. 
Since each coordinate $y_i$, for $i \in [1,2m]$, appears exactly once in each product, the product itself solves the heat equation.

\vspace{.1in}

\noindent
{\bf Step 3.  BC for $\Phi$.}
Consider the  cell faces 
$F_{i,2m}= \{\by: y_1< \ldots < y_i = y_{i+1} < \ldots < y_{2m}\}$ for $i \in [1,2m-1]$.
On $F_{2m,i}$ the Pfaffian 
reduces to a Pfaffian with $2n+2(m-1)$ rows and columns, namely where the row and column 
indexed by $(t,y_i)$ and $(t,y_{i+1})$ are removed. This can be seen since these two rows
and columns become identical except for the the entries $0,+1,-1$ where they cross. Then one may
 subtract row and column $(t,y_{i+1})$ from row and column $(t,y_i)$
 (using $\Pf(EBE^T) = \Pf(B) \det(E)$ for the corresponding elementary matrix) 
and expanding the Pfaffian along row $(t,y_i)$ leads to the Pfaffian of smaller size. Thus
\[
\Phi_{(t_i,z_i)}(t,\by)  =  \Phi_{(t_i,z_i)}(t,\by^{i,i+1})
\quad \mbox{for $\by \in F_{i,2m}, \; t > t_n$}
\]
where $\by^{i,i+1}=(y_1,\ldots,y_{i-1},y_{i+1},\ldots,y_{2m})$ (where $y^{i,i+1} = \emptyset $ if $m=2$). 

\vspace{.1in}

\noindent
{\bf Step 4.  $\Phi$ as an intensity.}
Fix smooth compactly supported $\phi_i:\R \to \R$, for $i \in [1,n]$, so that the 
distributions $\delta_{t_i} \times \phi_i$ are disjointly supported for $i \in [1,n]$. Consider the integral
\[
u(t,\by) = \int_{\R^n} \left( \prod_{i=1}^n \phi_i(z_i) dz_i \right)  \Phi_{(t_i,z_i)}(t,\by).
\]
The regularity of $\Phi_{(t_i,z_i)}$ above implies that
$u$ is a bounded $C^{1,2}((t_n,\infty) \times V_{2m}) \cap C([t_n,\infty) \times \overline{V}_{2m})$
solution to the heat equation with initial condition
\begin{equation} \label{ICu}
u(t_n,\by) =  \int_{\R^n} \left( \prod_{i=1}^n \phi_i(z_i) dz_i \right)  \Phi_{(t_i,z_i)}(t_n,\by)
\quad \mbox{for $\by \in V_{2m}$,}\\
\end{equation}
and boundary conditions
\begin{equation} \label{BCu}
u(t,\by)  =  \int_{\R^n} \left( \prod_{i=1}^n \phi_i(z_i) dz_i \right)  \Phi_{(t_i,z_i)}(t,\by^{i,i+1})
 \quad  \mbox{for $\by \in F_{i,2m}, \; t \geq t_n$.}
\end{equation}

\vspace{.1in}

\noindent
{\bf Step 5.  ABM mixed spin correlation and particle intensities.} Consider
\[
v(t,\by) = E \left[ \prod_{i=1}^n N_{t_i} (\phi_i) \; \prod_{j=1}^{2m} S_t(y_j) \right].
\]
The simple moment bounds $E[|N_t([a,b])|^k] < C(t,k) |b-a|^k $ from \cite{ourpaper} imply 
$v$ is a bounded function in $C([t_n,\infty) \times \overline{V}_{2m})$.
We claim also that $v \in C^{1,2}((t_n,\infty) \times V_{2m})$ and satisfies the heat
equation on  $(t_n,\infty) \times V_{2m}$. One way to see this is to apply the
time-duality from section 2.2 of \cite{ourpaper} to rewrite
\[
E \left[ \prod_{j=1}^{2m} S_t(y_j) | \sigma(N_s: s \leq t_n) \right]
=
E_{(y_1,\ldots,y_{2m})} \left[ (-1)^{\mu([\hat{X}_{t-t_n}^1, \hat{X}_{t-t_n}^2] \cup \ldots \cup
[\hat{X}_{t-t_n}^{2K-1}, \hat{X}_{t-t_n}^{2K}])} \right] |_{\mu = N_{t_n}}
\]
where the right hand side is the expectation over an annihilating system of Brownian motions
started from $(y_1,\ldots,y_{2m})$ which yields a set of $2K$ remaining particles at time $t-t_n$, 
positioned at $\hat{X}_{t-t_n}^1 < \ldots < \hat{X}_{t-t_n}^{2K}$ (here $K$ is random and possibly zero).
This satisfies the heat equation in $(t,\by)$ and can be used to show the same for $v$. 

Note that $v$ has initial condition
\begin{equation} \label{ICv}
v(t_n,\by) =  E \left[ \prod_{i=1}^n N_{t_i} (\phi_i) \, \prod_{j=1}^{2m} S_{t_n}(y_j) \right]
\quad \mbox{for $\by \in V_{2m}$,}\\
\end{equation}
and boundary conditions
\begin{equation} \label{BCv}
v(t,\by)  =  E \left[ \prod_{i=1}^n N_{t_i} (\phi_i) \!\!  \prod_{j=1,j \neq i,i+1}^{2m} \! S_t(y_j) \right]
 \quad  \mbox{for $\by \in F_{i,2m}, \; t \geq t_n$.}
\end{equation}

\vspace{.1in}

\noindent
{\bf Step 6. Induction.}
We aim to argue inductively that the initial and boundary conditions for $u$ and $v$ agree
(it is enough to consider boundary conditions only on the faces). 
This implies that $u=v$ and confirms that $\Phi_{(t_i,z_i)}$ is the desired intensity, completing the proof. 
In order to use induction we relabel the points $(t_i,z_i)$ as follows. Choose
$0<s_1 < \ldots < s_k$ and smooth compactly supported $(\phi_{i,i'}: i \in [1,k], i' \in [1,n_i])$, where $n_i \geq 1$,
so that the distributions $(\delta_{t_i} \times \phi_{i}: i \in[1,n])$ are precisely $(\delta_{s_i} \times \phi_{i,i'} :i \in [1,k], i' \in[1, n_i])$.  
We argue inductively first in $k \geq 0$, with a second inner induction on the number $2m$ of spin space points.
Note that when $k=0$, and the expectation has only spins, 
the theorem corresponds exactly to the Pfaffian found for $E[S_t(y_1) \ldots S_t(y_{2m})]$ in \cite{ourpaper}.

The initial condition (\ref{ICv}) for $v$ can be rewritten, noting that $s_k = t_n$, as
\begin{eqnarray*}
&& \hspace{-.3in} E \left[ \prod_{i=1}^k \prod_{i'=1}^{n_i} N_{s_i} (\phi_{i,i'}) \; \prod_{j=1}^{2m} S_{t_n}(y_j) \right] \\
& = &  \lim_{\hat{z}_l \downarrow z_l} E \left[ \prod_{i=1}^{k-1} \prod_{i'=1}^{n_i}
N_{s_i} (\phi_{i,i'}) \; \prod_{j=1}^{2m} S_{t_n}(y_j) 
\left(\prod_{l=1}^{n_k} \int \phi_{k,l}(z_l) S_{t_n}(z_l) S_{t_n}(\hat{z}_l) N_{t_n}(dz_l) \right) \right] \\
& = & 2^{-n_k}
\lim_{\hat{z}_l \downarrow z_l} E \left[ \prod_{i=1}^{k-1} \prod_{i'=1}^{n_i}
N_{s_i} (\phi_{i,i'}) \; \prod_{j=1}^{2m} S_{t_n}(y_j) 
\left(\prod_{l=1}^{n_k} \int \phi^{'}_{k,l}(z_l) S_{t_n}(z_l) S_{t_n}(\hat{z}_l) dz_l \right) \right] \\
& = & 2^{-n_k}
\lim_{\hat{z}_l \downarrow z_l}   \int_{R^{n_k}} \left( \prod_{l=1}^{n_k} \phi^{'}_{k,l}(z_l) \, dz_l \right) 
E \left[ \prod_{i=1}^{k-1} \prod_{i'=1}^{n_i} 
N_{s_i} (\phi_{i,i'})  \prod_{j=1}^{2m} S_{t_n}(y_j) \prod_{l=1}^{n_k}
S_{t_n}(z_l) S_{t_n}(\hat{z}_l) \right].
\end{eqnarray*} 
The distributional derivative used in the second equality, $(d/dz) S_t(z) = -2 S_t(z) N_t(dz)$, 
 holds almost surely and can be justified as in section 4.3 of \cite{ourpaper}.
The final expectation can be evaluated using the inductive hypothesis in $k$ in terms of the 
intensities $\Phi$. In brief, one can move the derivatives back from the test functions
$\phi_{k,l}$ onto these intensities, then carry out the limits $\hat{z}_l \downarrow z_l$, and one  
reaches the corresponding expression for the initial condition (\ref{ICu}) for $u$, 
confirming the initial conditions for $u$ and $v$ coincide. More carefully, 
note that the final expectation 
\begin{equation} \label{finalE}
E \left[ \prod_{i=1}^{k-1} \prod_{i'=1}^{n_i} 
N_{s_i} (\phi_{i,i'})  \prod_{j=1}^{2m} S_{t_n}(y_j) \prod_{l=1}^{n_k}
S_{t_n}(z_l) S_{t_n}(\hat{z}_l) \right]
\end{equation}
is a continuous function of the variables  $(z_1,\ldots,z_{n_k})$.
The integral over $\R^{n_k}$ can be broken into finitely many disjoint 
regions, according to the ordering of the points $(y_i:i \in [1,2m]) \cup (z_j,\hat{z}_j:j \in [1,n_k])$
on the real line. In each such region the inductive hypothesis gives a Pfaffian expression for (\ref{finalE}).
These expressions have continuous bounded derivatives in the variables $(z_j)$. 
Hence we may integrate by parts to take the derivatives off the test functions
$\phi^{'}_{k,l}$ and onto the Pfaffian expressions. By the boundedness of the derivatives
and the compact support of the test functions we may take the limits $\hat{z}_j \downarrow z_j$
inside the integral. This leaves a sum of Pfaffian expressions indexed over a smaller number of regions, 
namely the possible orderings of the points
$(y_i:i \in [1,2m]) \cup (z_j:j \in [1,n_k])$.
But in each of these regions the Pfaffian expression is identical, since 
switching two of the points $y_i$ and $z_j$, or $z_i$ and $z_j$, involves interchanging
two rows and columns leaving the Pfaffian unchanged.  Moreover, it is straightforward to check this 
single remaining Pfaffian is exactly the corresponding expression for the initial condition (\ref{ICu}) for $u$. 

The above argument shows that the statement of Theorem \ref{T2}
follows for $k$ and $m=0$ from the case $k-1$. 
An induction in $m$ implies immediately that the boundary condition (\ref{BCv}) for $v$ agrees with the boundary 
condition (\ref{BCu}) for $u$ and completes the double induction. 

\vspace{.1in}

\noindent
{\bf Acknowledgement.} We are grateful to Neil O'Connell for many illuminating discussions.

%%%%%%%%%%%%%%%%%%%%%%%%%%%%%%%%%%%%%%%%%%%%%%%%%%%%%%%

\end{document}